\documentclass[a4paper,11pt]{amsart}

\usepackage{amsmath,amsthm,amssymb,mathrsfs}
\usepackage{color}
\usepackage{xcolor}

\usepackage[colorlinks, linkcolor=blue,anchorcolor=Periwinkle,
    citecolor=blue,urlcolor=Emerald]{hyperref}
\usepackage[all]{xy}

\usepackage{pifont}
\usepackage{graphicx}
\usepackage{psfrag}

\usepackage{setspace}

\usepackage{tikz-cd}
\usepackage{bbm}

\usepackage{shorttoc}

\usepackage[normalem]{ulem}
\usepackage{tikz}
\usetikzlibrary{matrix,positioning,decorations.markings,arrows,decorations.pathmorphing,	
	backgrounds,fit,positioning,shapes.symbols,chains,shadings,fadings,calc}
\tikzset{->-/.style={decoration={  markings,  mark=at position #1 with
			{\arrow{>}}},postaction={decorate}}}
\tikzset{-<-/.style={decoration={  markings,  mark=at position #1 with
			{\arrow{<}}},postaction={decorate}}}

\numberwithin{equation}{section}

\newcommand{\ie}{{\em i.e.}\ }
\newcommand{\cf}{{\em cf.}\ }

\numberwithin{equation}{section}
\newtheorem{theorem}{Theorem}[section]
\newtheorem*{theorem*}{Theorem}

\newtheorem{lemma}[theorem]{Lemma}
\newtheorem{notation}[theorem]{Notations}
\newtheorem{proposition}[theorem]{Proposition}
\newtheorem{corollary}[theorem]{Corollary}

\newtheorem*{conjecture*}{Conjecture}

\newtheorem{remark}[theorem]{Remark}
\newtheorem{definition}[theorem]{Definition}

\newcommand{\Hom}{\mathrm{Hom}}

\newcommand{\Z}{\mathbb{Z}}

\newcommand{\id}{\mathbf{1}}

\newcommand{\supp}{\mathrm{supp}}

\renewcommand{\phi}{\varphi}

\renewcommand{\hat}[1]{\widehat{#1}}

\renewcommand{\tilde}[1]{\widetilde{#1}}

\newcommand{\diag}{\text{diag}}

\makeatletter
\DeclareRobustCommand{\cev}[1]{%
  {\mathpalette\do@cev{#1}}%
}
\newcommand{\do@cev}[2]{%
  \vbox{\offinterlineskip
    \sbox\z@{$\m@th#1 x$}%
    \ialign{##\cr
      \hidewidth\reflectbox{$\m@th#1\vec{}\mkern4mu$}\hidewidth\cr
      \noalign{\kern-\ht\z@}
      $\m@th#1#2$\cr
    }%
  }%
}
\makeatother

\title[$f$-vectors and $F$-invariant in generalized cluster algebras]{$f$-vectors and $F$-invariant in generalized cluster algebras}

\author[C. Fu]{Changjian Fu}
\address{Department of Mathematics, Sichuan University, Chengdu 610064, P.R.China}
\email{changjianfu@scu.edu.cn}
\author[H. Ye]{Huihui Ye}
\address{Department of Mathematics, Sichuan University, Chengdu 610064, P.R.China}
\email{YeHuiHuimath@outlook.com}

\subjclass[2020]{}

\begin{document}
\begin{abstract}
We establish certain fundamental properties of $f$-vectors and $F$-matrices for generalized cluster algebras, including the initial and final seed mutation formulas, the compatibility property and the symmetry property. Along the way, we also generalize the construction of $F$-invariant for generalized cluster algebras without assuming positivity and prove certain basic properties.

\end{abstract}
\maketitle
\tableofcontents
\section{Introduction}
Fomin and Zelevinsky invented cluster algebras in \cite{ fomin_zelevinsky_2002} to provide a combinatorial framework for studying total positivity in algebraic groups and canonical bases of quantum groups. Since then, cluster algebras have been found to have deep connections with many other areas of mathematics and physics, such as discrete dynamical systems, non-commutative algebraic geometry, string theory and quiver representation theory, \cf \cite{Keller2012} and the references therein.

A cluster algebra is a commutative algebra that possesses a unique set of generators known as cluster variables. These generators are gathered into
overlapping sets of fixed finite cardinality, called clusters, which are defined recursively from an initial one via a mutation operation. The exchange matrix determines the mutations of clusters in different directions. 
A compatibility degree of cluster algebra is a function on the set of pairs of cluster variables satisfying various properties. This function was first introduced by Fomin and Zelevinsky \cite{fomin_zelevinsky_2003b} for generalized associahedra associated with finite root systems in their study of Zamolodchikov’s periodicity for $Y$-systems, which are a special kind of cluster complexes of cluster algebras. It has played a key role in the classification of cluster algebras of finite type. As a generalization of the classical compatibility degree of Fomin and Zelevinsky, Cao and Li \cite{Cao_Li20} introduced the $d$-compatibility degree for any cluster algebras by using $d$-vectors. However, the $d$-compatibility degree does not preserve many properties of the classical one. 
The first named author and Gyoda \cite{Fu_Gyoda} introduced another generalization, called $f$-compatibility degree, for any cluster algebras by using $f$-vectors, and established various interesting properties for it.
As a further generalization,
Cao \cite{cao23} introduced the $F$-invariant in cluster algebras using tropicalization, which is an analog of $E$-invariant introduced by Derksen-Weyman-Zelevinsky \cite{Derksen_Weyman_Jerzy_Zelevinsky10} in the additive categorification of cluster algebras using decorated representations of the quiver with potentials
and the $\sigma$-invariant introduced by Kang-Kashiwara-Kim-Oh \cite{Kang_Kashiwara_Kim_Oh18} in the monoidal categorification of (quantum) cluster algebras, using representations of quiver Hecke algebras.

 In their study of Teichm\"{u}ller space of Riemann surface with orbifold points, Chekhov and Shapiro \cite{Chekhov_Shapiro} introduced {\it generalized cluster algebras}, which is a significant generalization of classical cluster algebras. Chekhov and Shapiro \cite{Chekhov_Shapiro} prove the Laurent phenomenon for generalized cluster algebras, while Nakanishi \cite{Nakanishi14} demonstrated that these algebras share fundamental structural parallels with classical cluster algebras. Specifically, Nakanishi introduced \textit{F}-polynomials and established the separation formulas for generalized cluster algebras. Thanks to separation formulas, all cluster variables can describled by the {\it $C$-matrices}, the {\it $G$-matrices} and the {\it $F$-polynomials}, where the {$C$-} and {$G$-matrices} are the tropical part and the {\it $F$-polynomials} are nontropical part. It's worth mentioning that the structure of generalized cluster algebras also appears in many other branches of math, such as the representation theory of quantum affine algebra \cite{Gleitz2015}, WKB analysis \cite{Iwaki-Nakanishi2016} and representation theory of finite dimensional algebras \cite{LFM18,LF-Velasco2018}. 

In this paper, we define $f$-vectors/$F$-matrices for generalized cluster algebras and obtain the initial and final seed mutations of $F$-matrices.  It is proved that $f$-vectors/ $F$-matrices  satisfy  the symmetry property and the compatibility property (\cf Theorem \ref{relation of f-compatible and cluster}). 
 We also introduce the $F$-invariant for upper generalized cluster algebras without assuming positivity.  The $F$-invariant $(u\parallel u')_F$ vanishes (\ie, $(u\parallel u')_F=0$) for any cluster monomials $u$ and $u'$ in the same cluster.
Moreover, a good element $u$ is a cluster monomial if and only if there exists a vertex $t\in\mathbb{T}_n$ such that $(x_{i;t}\parallel u)_F=0$ for any $i\in [1,n]$. 
Unfortunately, the inverse proposition is not right. Because there is a good element $u$ in the upper cluster algebra $\mathcal{U}$ such that $(u\parallel u)_F=0$, but $u$ is not a cluster monomial, \cf \cite[Example 4.22]{cao23}. 
In the categorification of generalized cluster algebras arising from surfaces with orbifold points of order 3, Daniel Labardini-Fragoso and Lang Mou \cite{LFM18} gives bijection from $\tau$-rigid pairs and cluster monomials of the generalized cluster algebras using the Caldero-Chapoton map. It would be interested to compare the $E$-invariants and $F$-invariants in this setting, which will be discussed in the future.

This paper is organized as follows. In Section \ref{preliminaries}, we recall some basic definitions, notation, and results for generalized cluster algebras. In Section \ref{f-vectors}, after introducing the $f$-vectors/$F$-matrices for generalized cluster algebras, we prove the fundamental properties of $f$-vectors (Theorem \ref{relation of f-compatible and cluster}),
 and establish the initial-seed mutation formulas of $F$-matrices (Proposition \ref{prop:initial-seed-mutation}). 
In Section \ref{F-invariant in generalized cluster algebras}, we define the $F$-invariant for generalized cluster algebras without assuming positivity (compare to \cite{cao23}) and establish certain fundamental results for it (\cf Corollary \ref{criterion for the same cluster of cluster variables and good elements} and Proposition \ref{uu' is cluster monomial=(u,u')_F=0}).


\section{Preliminaries}\label{preliminaries}
\subsection{Generalized cluster algebras}
In this section, we recall some basics of generalized cluster algebras \cite{Chekhov_Shapiro,Nakanishi14}.
Fix positive integers $m\ge n$ and an n-tuple $\mathbf{r}=(r_{1},\dots,r_{n})$ of positive integers. We denote by $[1,n]=\{1,2.\dots,n\}$ and $R=\operatorname{diag}\{r_1,\dots, r_n\}$.
Let  $\mathbf{z}=(z_{i,s_i})_{i=1,2,\dots,n;s_i=1,2,\dots,r_i-1}$ with $z_{i,s_i}=z_{i,r_{i}-s_i}$ be formal variables and $\mathcal{F}=\mathbb{Q}(\mathbf{z},x_1,\dots, x_m)$ the rational function field in $\mathbf{z}$ and $x_1,\dots, x_m$.
Denote by \[\mathcal{F}_{>0}:=\mathbb{Q}_{sf}(z_{i,s_i};x_1,x_2,\dots,x_m|1\le i\le n,1\le s_i\le r_{i}-1)\] the universal semi-field in $\mathbf{z}$ and $x_{1},\dots,x_m$. 
In this paper, we fix the semi-field $\mathbb{P}=\operatorname{Trop}(\{x_{n+1},\dots,x_{m}\},\mathbf{z})$ the tropical semi-field generated by the elements $\{x_{n+1},\dots,x_{m}\}$ and $\mathbf{z}$, which is the multiplicative abelian group with tropical sum $\oplus$ defined by
\begin{equation*}
(\prod\limits_{i}x_i^{a_i}\prod\limits_{i,s}z_{i,s}^{a_{i,s}})\oplus(\prod\limits_{i}x_i^{b_i}\prod\limits_{i,s}z_{i,s}^{b_{i,s}})=(\prod\limits_{i}x_i^{\min\{a_i,b_i\}}\prod\limits_{i,s}z_{i,s}^{\min\{a_{i,s},b_{i,s}\}}),
\end{equation*}
where $a_i,a_{i,s},b_i,b_{i,s}\in \Z$. Let $\mathbb{ZP}$ be the group ring of $\mathbb{P}$ and $\mathbb{QP}$ be the skew field of $\mathbb{ZP}$. 
A {\it compatible pair} $(\tilde{B}, \Lambda)$ consists of an integer $m\times n$-matrix $\tilde{B}$ and a skew-symmetric integer $m\times m$-matrix $\Lambda$ such that 
\[
\tilde{B}^T\Lambda=[D\ 0],
\]
where $D=\operatorname{diag}(d_1,\dots, d_n)$ is a diagonal $n\times n$ matrix whose diagonal coefficients are positive integers. It is straightforward to verify that the principal part $B$ (\ie the submatrix formed by the first $n$ rows) of $\tilde{B}$ is skew-symmetrizable and $D$ is a skew-symmetrizer of $B$.
\begin{remark}
   It is convenient  to extend $\mathbf{z}$ to $\mathbf{z}=(z_{i,s_i})_{i=1,\dots, n;s_i=0,\dots, r_i}$ (by abuse notation, we also denote it by $\mathbf{z}$) by setting $z_{i,0}=z_{i,r_i}=1$. The pair $(\mathbf{r},\mathbf{z})$ is referred to as the {\it mutation data}.  Throughout this subsection, we fix the {\it mutation data} $(\mathbf{r},\mathbf{z})$. 
\end{remark}

\begin{definition}
      A {\it labeled $(\mathbf{r},\mathbf{z})$-seed} is a pair $(\mathbf{x},\tilde{B})$ such that. 
    \begin{itemize}
        \item $\tilde{B}=(b_{ij})$ is a $m\times n$ integer matrix, of which the principal part $B$ is skew-symmetrizable and $D$ is a skew-symmetrizer of $B$;
        \item $\mathbf{x}=(x_1,\dots,x_m)$ is an m-tuple of algebraic independent elements of $\mathcal{F}$ over the rational function field $\mathbb{Q}(\mathbf{z})$ of $\mathbf{z}$.
    \end{itemize}
We say that $\mathbf{x}$ is an {\it $(\mathbf{r},\mathbf{z})$-cluster} and refer to $x_i$ and $\tilde{B}$ as the {\it cluster variables}
and the {\it exchange matrix}, respectively.
\end{definition}
\begin{definition}
     The $(\mathbf{r},\mathbf{z})$-$Y$-seed of rank n in $\mathcal{F}$ is a pair $(\mathbf{y},\hat{B})$, where
     \begin{itemize}
         \item $\mathbf{y}=(y_1,\dots,y_{m})$ is a freely generating set of $\mathcal{F}$ over $\mathbb{Q}(\mathbf{z})$,
         \item $\hat{B}=(B|Q)=(\hat{b}_{ij})$ is an $n\times m$ integer matrix such that $B$ is a skew-symmetrizable matrix.
     \end{itemize}
 \end{definition}
Let $k\in [1,n]$.
We define $E_{k,\varepsilon}^{\tilde{B}R}$ as the $m\times m$-matrix which differs from the identity matrix only in its $k$-th column whose coefficients are given by
    \begin{align*}
       (E_{k,\varepsilon}^{\tilde{B}R})_{ik}=
   \begin{cases}
         -1& \text{if } i=k;\\
           [-\varepsilon b_{ik}r_k]_+& \text{if } i\neq k,
       \end{cases}  
     \end{align*}
   where $\varepsilon\in\{\pm1\}$.
 Denote by $F_{k,\varepsilon}^{R\tilde{B}}$ the $n\times n$-matrix that differs from the identity matrix only in its $k$-th row, with coefficients given by
    \begin{align*}
        (F^{R\tilde{B}}_{k,\varepsilon})_{ki}=\begin{cases}
            -1& if\ i=k;\\
           [\varepsilon r_{k}b_{ki}]_+& if\ i\neq k,
       \end{cases}
    \end{align*}
   where $\varepsilon\in\{\pm 1\}$.

The {\it mutation $\mu_k$ in direction $k$}  transforms the compatible pair $(\tilde{B},\Lambda)$  into $\mu_{k}(\tilde{B},\Lambda):=(\tilde{B}',\Lambda')$, where
\[
\tilde{B}'=E_{k,\varepsilon}^{\tilde{B}R}\tilde{B}F^{R\tilde{B}}_{k,\varepsilon},\ \Lambda'=(E_{k,\varepsilon}^{\tilde{B}R})^T\Lambda E_{k,\varepsilon}^{\tilde{B}R}.
\]
 In fact, the mutation of a compatible pair is an involution.

Now we recall the $(\mathbf{r},\mathbf{z})$-mutation in generalized cluster algebras.
\begin{definition}\label{Definition of mutation}
    For any $(\mathbf{r},\mathbf{z})$-seed ($\mathbf{x}$,$\tilde{B}$) and $k\in[1,n]$, the $(\mathbf{r},\mathbf{z})$-mutation of $(\mathbf{x},\tilde{B})$ in direction $k$ is a new $(\mathbf{r},\mathbf{z})$-seed $\mu_k(\mathbf{x},\tilde{B}):=(\mathbf{x}',\tilde{B}') $ defined by the following rule:
\begin{align}
x_i'&=
\begin{cases}
        x_i & \text{if $i\neq k$;} \\
         x_k^{-1}(\prod\limits_{j=1}^{m}x_j^{[-\varepsilon b_{jk}]_+})^{r_k}\sum\limits_{s=0}^{r_k}z_{k,s}\hat{y}_k^{\varepsilon s}
         & \text{if $i=k$};
\end{cases}\label{E4}\\
b_{ij}'&=
 \begin{cases}
-b_{ij} &\text{if $i=k$ or $j=k$};\\
b_{ij}+r_k([-\varepsilon b_{ik}]_+b_{kj}+b_{ik}[\varepsilon b_{kj}]_+) &\text{else}.
\end{cases}\label{E6}
\end{align} 
where $\varepsilon\in\{\pm 1\}$ and $\hat{y}_{k}=\prod\limits_{j=1}^{m}x_{j}^{b_{jk}}$.
\end{definition}

\begin{remark}\label{R1}
\begin{enumerate}
    \item The mutation formulas (\ref{E4}) and (\ref{E6})  are independent of the choice of $\varepsilon$, and $\mu_k$ is an involution;
\item $x_{1},\dots,x_{n}$ is called unfrozen cluster variables and $x_{n+1},\dots,x_{m}$ is called frozen variables or coefficients;
    \item If $\mathbf{r}=(1,\dots,1)$, then the mutation formulas (\ref{E4}) and (\ref{E6}) reduce to the mutation formulas of cluster algebras.
\end{enumerate}


\end{remark}

\begin{definition}
    The mutation of $(\mathbf{r},\mathbf{z})$-$Y$-seed $(\mathbf{y},\hat{B})$ in direction $k\in[1,n]$ is the pair $(\mathbf{y}',\hat{B}'):=\mu_{k}(\mathbf{y},B)$ given as follows:
   
\begin{align}
            y_{i}'&=\begin{cases}
                y_{k}^{-1} & if\quad i=k,\\
                y_{i}(y_{k}^{[\varepsilon \hat{b}_{ki}]_{+}})^{r_{k}}(\sum\limits_{s=0}^{r_{k}}z_{k,s}y_{k}^{\varepsilon s})^{-\hat{b}_{ki}} & if\quad i\neq k,
            \end{cases}\\
            \hat{b}_{ij}'&=\begin{cases}
                   - \hat{b}_{ij}  &if\quad i=k\ \ or\ \ j=k,\\
                \hat{b}_{ij}+r_{k}([-\varepsilon\hat{b}_{ik}]_{+}\hat{b}_{kj}+\hat{b}_{ik}[\varepsilon \hat{b}_{kj}]_{+}) & if \quad i,j\neq k.
            \end{cases}
        \end{align}     
The variables $y_1,\dots,y_m$ are called $y$-variables of $(\mathbf{y},\hat{B})$.
\end{definition}

The mutation $(\mathbf{y}',\hat{B}')$ of $(\mathbf{r},\mathbf{z})$-$Y$-seed in direction $k$ is also an $(\mathbf{r},\mathbf{z})$-$Y$-seed and $\mu_{k}$ is an involution.

Let $\mathbb{T}_n$ be the $n$-regular tree. An $(\mathbf{r},\mathbf{z})$-cluster pattern $\mathbf{\Sigma}=\{\Sigma_t:=(\mathbf{x}_t,\tilde{B}_t)\}_{t\in \mathbb{T}_n}$ of rank $n$ is a collection of $(\mathbf{r},\mathbf{z})$-seeds $(\mathbf{x}_t,\tilde{B}_t)$ in $\mathcal{F}$ such that $\Sigma_{t'}=\mu_k(\Sigma_t)$ whenever $t,t'$ are $k$-adjacent in $\mathbb{T}_n$. Similarly, one can define the $(\mathbf{r},\mathbf{z})$-$Y$-pattern consisting of $(\mathbf{r},\mathbf{z})$-$Y$-seed. 
For an $(\mathbf{r},\mathbf{z})$-cluster pattern $\mathbf{\Sigma}=\{\Sigma_t:=(\mathbf{x}_t,\tilde{B}_t)\}_{t\in \mathbb{T}_n}$, we always denote by $\mathbf{x}_t=(x_{1;t},\dots, x_{m;t})$, $\tilde{B}_t=(\mathbf{b}_{j;t})=(b_{ij;t})$ and by $B_t$ the principal part of $\tilde{B}_t$. We also refer to  $\{\tilde{B}_t\mid t\in \mathbb{T}_n\}$ as the $(\mathbf{r},\mathbf{z})$-${B}$-pattern. We fix a vertex $t_0$ as {\it root vertex}, i.e., an $(\mathbf{r},\mathbf{z})$-cluster pattern is uniquely determined by assigning $(\mathbf{x}_{t_0}, \tilde{B}_{t_0})$ to $t_0$, and refer to $\Sigma_{t_0}$ as the {\it initial seed}.

Following \cite[Definition 2.6]{cao23}, we recall the definition of a Langland-Poisson triple.
\begin{definition}
    Let  $\mathcal{S}_{X}=\{(\mathbf{x}_{t},\tilde{B}_{t})\}_{t\in \mathbb{T}_n}$ be an $(\mathbf{r},\mathbf{z})$-cluster pattern, $\mathcal{S}_{Y}=\{(\mathbf{y}_{t},\hat{B}_{t})\}_{t\in\mathbb{T}_n}$ an $(\mathbf{r},\mathbf{z})$-$Y$-pattern and $\Lambda=\{\Lambda_{t}|t\in\mathbb{T}_{n}\}$ be a collection of skew-symmetric matrices indexed by the vertices in $\mathbb{T}_{n}$.
    \begin{enumerate}
        \item The pair $(\mathcal{S}_{X},\mathcal{S}_{Y})$ is called a {\it Langland dual pair} if $\hat{B}_{t_{0}}=-\tilde{B}_{t_{0}}^{\top}$ holds. 
        \item The triple $(\mathcal{S}_{X},\mathcal{S}_{Y},\Lambda)$ is called a {\it Langland-Poisson triple} if it satisfies that 
\begin{itemize}
    \item $(\mathcal{S}_{X},\mathcal{S}_{Y})$ is {\it Langland dual pair};
    \item $\{(\tilde{B}_{t},\Lambda_{t})|t\in\mathbb{T}_{n}\}$ forms a collection of compatible pairs and $(\tilde{B}_{t'},\Lambda_{t'})=\mu_{k}(\tilde{B}_{t},\Lambda_{t})$ whenever $t$ and $t'$ are $k$-adjacent in $\mathbb{T}_n$.
\end{itemize}
\end{enumerate}
\end{definition}
\begin{definition}
\label{Definition of generalized cluster algebras}
Let $\mathbf{\Sigma}=\{\Sigma_t:=(\mathbf{x}_t,\tilde{B}_t)\}_{t\in \mathbb{T}_n}$ be an $(\mathbf{r},\mathbf{z})$-cluster pattern of rank $n$.
     The {\it generalized cluster algebra} $\mathcal{A}:=\mathcal{A}(\mathbf{\Sigma})$ associated to  $\mathbf{\Sigma}$  is the $\Z \mathbb{P}$-subalgebra of $\mathcal{F}$ generated by $\mathcal{X}(\mathbf{\Sigma}):=\bigcup_{t\in \mathbb{T}_n}\{x_{1;t},\dots, x_{n;t}\}$.
     The upper generalized cluster algebra $\mathcal{U}$ of $\mathbf{\Sigma}$ is the $\mathbb{ZP}$-subalgebra of $\mathcal{F}$ given by
     \[\mathcal{U}:=\bigcap_{t\in\mathbb{T}_{n}}\mathcal{L}(t),\]
    where $\mathcal{L}(t)=\mathbb{Z}\mathbb{P}[x_{1;t}^{\pm},\cdots,x_{m;t}^{\pm}]$.
\end{definition}
We also refer to generalized cluster algebras and upper generalized  cluster algebras as the {\it $(\mathbf{r},\mathbf{z})$-cluster algebras} and the {\it $(\mathbf{r},\mathbf{z})$-upper cluster algebras}, respectively.

The Laurent phenomenon still holds for generalized cluster algebras.
\begin{proposition}\cite[Theorem 2.5]{Chekhov_Shapiro}\label{Laurent phenomenon}
     Each cluster variable $x_{i;t}$ could be expressed as a Laurent polynomial in $\mathbf{x}$ with coefficients in $\mathbb{Z}\mathbb{P}$.
\end{proposition}

\subsection{$C$- and $G$-matrices}
We fix an $(\mathbf{r},\mathbf{z})$-cluster pattern $\mathbf{\Sigma}=\{\Sigma_t=(\mathbf{x}_t,\tilde{B}_t)\}_{t\in \mathbb{T}_n}$ with root vertex $t_0\in \mathbb{T}_n$. For each vertex $t\in \mathbb{T}_n$, we introduce two integer matrices $C_t=(\mathbf{c}_{1;t},\dots,\mathbf{c}_{n;t} )=(c_{ij;t})_{i,j=1}^n$ and $G_t=(\mathbf{g}_{1;t},\dots, \mathbf{g}_{n;t})=(g_{ij;t})_{i,j=1}^n$, called {\it $C$-matrix} and {\it $G$-matrix} respectively, by the following recursions:
\begin{itemize}
    \item  $C_{t_0}=G_{t_0}=I_n$;
    \item If $\xymatrix{t\ar@{-}[r]^k&t'}\in \mathbb{T}_n$, then
    \begin{align}
 \label{E9}
         c_{ij:t'}&=
        \begin{cases}
            -c_{ij;t} & \text{if $j=k$};\\
        c_{ij;t}+r_k(c_{ik;t}[\varepsilon b_{kj;t}]_++[-\varepsilon c_{ik;t}]_+b_{kj;t})& \text{if $ j\neq k$};
\end{cases}\\ 
\label{E10}
        \mathbf{g}_{i;t'}&=
        \begin{cases}
            \mathbf{g}_{i;t} & \text{if $i\neq k$};\\
-\mathbf{g}_{k;t}+r_k(\sum\limits_{j=1}^{n}[-\varepsilon b_{jk;t}]_+\mathbf{g}_{j;t}-\sum\limits_{j=1}^{n}[-\varepsilon c_{jk;t}]_+\mathbf{b}_{j;t_0}) & \text{if $i=k$}.
        \end{cases}
    \end{align}
    \end{itemize}
The recurrence formulas (\ref{E9}) and (\ref{E10}) are independent of the choice of the sign $\varepsilon\in \{\pm 1\}$.  The column vectors of $C_t$ and $G_t$ are called {\it $c$-vectors} and {\it $g$-vectors} of $\mathbf{\Sigma}$ respectively.
We refer to $\{C_t\}_{t\in \mathbb{T}_n}$ and $\{G_t\}_{t\in\mathbb{T}_n}$ as the {\it $(\mathbf{r},\mathbf{z})$-$C$-pattern} and {\it $(\mathbf{r},\mathbf{z})$-$G$-pattern} of $\mathbf{\Sigma}$ with respect to $t_0$. 
\begin{remark}
    Since $\{C_t\}_{t\in \mathbb{T}_n}$ and $\{G_t\}_{t\in\mathbb{T}_n}$ only depend on $B_{t_0}, \mathbf{r}$, and $t_0\in \mathbb{T}_n$, we also call  $\{C_t\}_{t\in \mathbb{T}_n}$ (resp. $\{G_t\}_{t\in\mathbb{T}_n}$) the $(\mathbf{r},\mathbf{z})$-$C$-pattern (resp. $(\mathbf{r},\mathbf{z})$-$G$-pattern) associated with the $(\mathbf{r},\mathbf{z})$-$B$-pattern $\{B_t\mid t\in \mathbb{T}_n\}$ with respect to $t_0$. We also write $C_t^{B_{t_0};t_0}$ and $G_t^{B_{t_0};t_0}$ for $C_t$ and $G_t$, respectively, to emphasize the exchange matrix $B_{t_0}$ at $t_0$. 

    When $\mathbf{r}=(1,\dots, 1)$, the $(\mathbf{r},\mathbf{z})$-$B$-pattern becomes an ordinary $B$-pattern, and the $(\mathbf{r},\mathbf{z})$-$C$-pattern  and $(\mathbf{r},\mathbf{z})$-$G$-pattern degenerate to ordinary $C$-pattern and $G$-pattern respectively.
\end{remark}

\begin{proposition}\cite{Nakanishi14}\label{relation_of_C-matriice} The following relations holds:
    \begin{align}
        C_{t}&=^{L}C_{t}=R(^{R}C_{t})R^{-1},\\
        G_{t}&=^{R}G_{t}=R^{-1}(^{L}G_{t})R,
    \end{align}
    where $\{^LC_t\}_{t\in \mathbb{T}_{n}}$ and $\{^LG_t\}_{t\in \mathbb{T}_{n}}$ are the ordinary $C$-pattern and the ordinary $G$-pattern associated with ordinary $B$-pattern $\{ RB_{t}|t\in\mathbb{T}_{n}\}$ with respect to $t_0$, while $\{^RC_t\}_{t\in \mathbb{T}_{n}}$ and $\{^RG_t\}_{t\in \mathbb{T}_{n}}$ are $C$-pattern and $G$-pattern associated with $B$-pattern $\{ B_{t}R|t\in\mathbb{T}_{n}\}$ with respect to $t_0$.
\end{proposition}

By the sign-coherence of $c$-vectors \cite[Corollary 5.5]{GHKK2018}, (\ref{E10}) can be rewritten as
\begin{align}\label{E15}
\mathbf{g}_{i;t'}=\begin{cases}\mathbf{g}_{i;t}& \text{if $i\neq k$;}\\
    -\mathbf{g}_{k;t}+r_k(\sum\limits_{j=1}^{n}[-\varepsilon_{k;t} b_{jk;t}]_+\mathbf{g}_{j;t})& \text{if $i=k$},
    \end{cases}
\end{align}
where $\varepsilon_{k;t}$ is the common sign of components of the $c$-vector $\mathbf{c}_{k;t}$.

For each vertex $t\in \mathbb{T}_n$, we assign an $m\times m$-integer matrix $\tilde{G}_t=(\tilde{\mathbf{g}}_{1;t},\dots, \tilde{\mathbf{g}}_{m;t})$ to $t$ by the following recursion:
\begin{itemize}
    \item $\tilde{G}_{t_0}=I_m$;
    \item If $\xymatrix{t\ar@{-}[r]^k&t'}\in \mathbb{T}_n$, then
     \begin{align}\label{eq:g-formula}
         \mathbf{\tilde{g}}_{i;t'}=
        \begin{cases}
            \mathbf{\tilde{g}}_{i;t} & \text{ if } i\neq k;\\
            -\mathbf{\tilde{g}}_{k;t}+r_{k}(\sum\limits_{j=1}^{m}[-b_{jk;t}]_+\mathbf{\tilde{g}}_{j;t}-\sum\limits_{j=1}^{n}[-  c_{jk;t}]_+\mathbf{b}_{j;t_0}) & \text{ if }i=k.
        \end{cases}
    \end{align}

\end{itemize}
We call $\tilde{G}_t$ the {\it extended $G$-matrix} of $\mathbf{\Sigma}$ and its column vectors the {\it extended $g$-vectors}. Obviously, each extended $G$-matrix takes the form
\begin{align*}
    \tilde{G}_t=\begin{bmatrix}
        G_{t}&0\\
        \star & I_{m-n}
    \end{bmatrix}.
\end{align*}
\section{\texorpdfstring{F-matrices and their mutation formulas}{f-vectors}}\label{f-vectors}
In this section, we introduce $f$-vectors and $F$-matrices for generalized cluster algebras and investigate its properties. In particular, we show the recurrence formula and initial-seed mutation formula for $F$-matrices, and establish the fundamental properties for $f$-vectors.

Throughout this section, we fix an $(\mathbf{r},\mathbf{z})$-cluster pattern $\mathbf{\Sigma}=\{\Sigma_t=(\mathbf{x}_t,\tilde{B}_t)\}_{t\in \mathbb{T}_n}$ with root vertex $t_0\in \mathbb{T}_n$.
\subsection{$F$-matrices and their mutation formula}

\begin{definition}
    \label{Definition of F-polynomial}
For each $t\in \mathbb{T}_n$, the $i$th $F$-polynomial $F_{i;t}(\mathbf{\hat{y}}_{t_{0}},\mathbf{z})$ at $t$ is defined by the following recursion:
\begin{itemize}
    \item  $F_{i;t_0}(\mathbf{\hat{y}}_{t_0},\mathbf{z})=1$ for $i\in[1,n]$;
    \item  If $\xymatrix{t\ar@{-}[r]^k&t'}\in \mathbb{T}_n$, then  \begin{align*}
        F_{i;t'}(\mathbf{\hat{y}}_{t_0},\mathbf{z})&=\begin{cases}
            F_{k;t}(\mathbf{\hat{y}}_{t_0},\mathbf{z})^{-1}M_{k;t} & if\quad i=k,\\
            F_{i;t}(\mathbf{\hat{y}}_{t_0},\mathbf{z}) & if\quad i\neq k,
        \end{cases}
    \end{align*}
    where $M_{k;t}=(\prod\limits_{j=1}^{n}\hat{y}_{j;t_0}^{[-\varepsilon c_{jk;t}]_{+}}F_{j;t}(\mathbf{\hat{y}}_{t_0},\mathbf{z})^{[-\varepsilon b_{jk;t}]_{+}})^{r_{k}}\sum\limits_{s=0}^{r_{k}}z_{k,s}(\prod\limits_{j=1}^{n}\hat{y}_{j;t_0}^{\varepsilon c_{jk;t}}F_{j;t}(\mathbf{\hat{y}}_{t_0},\mathbf{z})^{\varepsilon b_{jk;t}})^{s}$ and $\{C_t=(c_{ij;t})\mid t\in \mathbb{T}_n\}$ is the {\it $(\mathbf{r},\mathbf{z})$-$C$-pattern} of $\mathbf{\Sigma}$.
\end{itemize}  
\end{definition}

According to \cite{Nakanishi14}, $F_{i;t}(\hat{\mathbf{y}}_{t_0},\mathbf{z})\in \mathbb{Z}[\hat{\mathbf{y}}_{t_0},\mathbf{z}]$ and the followings are established:
\begin{itemize}
    \item $F_{i;t}(\mathbf{\hat{y}}_{t_0},\mathbf{z})$ has a constant term 1;
    \item $F_{i;t}(\mathbf{\hat{y}}_{t_0},\mathbf{z})$ has a unique monomial $\mathbf{\hat{y}}_{t_0}^{\mathbf{f}^{t_0}_{i;t}}$ as a polynomial in $\mathbf{\hat{y}}_{t_0}$  such that each monomial $\mathbf{\hat{y}}_{t_0}^{a}$ with nonzero coefficient in $F_{i;t}(\mathbf{\hat{y}}_{t_0},\mathbf{z})$ divides $\mathbf{\hat{y}}_{t_0}^{\mathbf{f}^{t_0}_{i;t}}$, \ie, $F_{i;t}(\mathbf{\hat{y}}_{t_0},\mathbf{z})$ has a ``maximal degree" monomial. Furthermore, the coefficient of monomial $\mathbf{\hat{y}}_{t_0}^{\mathbf{f}^{t_0}_{i;t}}$ appears in $F_{i;t}(\mathbf{\hat{y}}_{t_0},\mathbf{z})$ is $1$.
\end{itemize}
\begin{definition}
The integer vector $\mathbf{f}_{i;t}^{t_0}$ is called the $i$th {\it f-vector} of $\mathbf{\Sigma}$ at $t$ with respect to the root vertex $t_0$ and the matrix $F_t^{t_0}:=(\mathbf{f}_{1;t}^{t_0},\dots,\mathbf{f}_{n;t}^{t_0})$ is called the {\it F-matrix} of $\mathbf{\Sigma}$ at $t$ with respect to $t_0$. The collection $\mathbf{F}=\{F_{t}^{t_0}\mid t\in \mathbb{T}_n\}$ is called the $F$-pattern associated to $\mathbf{\Sigma}$.
\end{definition}
\begin{remark}
    We may fix any vertex $s\in \mathbb{T}_n$ as a root vertex of $\mathbf{\Sigma}$ and the $f$-vector $\mathbf{f}_{i;t}^s$ depend on the choice of $s$.
\end{remark}
It is clear that $F$-polynomials $F_{i;t}(\hat{\mathbf{y}}_{t_0},\mathbf{z})$ and the $F$-pattern $\{F_t^{t_0}\mid t\in \mathbb{T}_n\}$ are uniquely determined by $B_{t_0}$ and the mutation data $(\mathbf{r},\mathbf{z})$. We also write $F_t^{B_{t_0};t_0}$ for $F_t^{t_0}$ to emphasize the exchange matrix $B_{t_0}$ at $t_0$.

    
Since the $(\mathbf{r},\mathbf{z})$-cluster algebra does not depend on the sign of the initial exchange matrix, the results in \cite[Theorem 2.8]{Fujiwara_Gyoda19} extend to generalized cluster algebras. 
\begin{theorem}\label{relation of C,G and F matrices}
We have the following relations:
    \begin{align}
        C_{t}^{-B_{t_0};t_0}&=C_{t}^{B_{t_0};t_0}+F_{t}^{B_{t_0};t_0}B_t,\label{eqn:c-matrix-duality}\\
        G_{t}^{-B_{t_0};t_0}&=G_{t}^{B_{t_0};t_0}+B_{t_0}F_{t}^{B_{t_0};t_0},\\
        F_{t}^{-B_{t_0};t_0}&=F_{t}^{B_{t_0};t_0.}
    \end{align}
    
\end{theorem}


Let $J_n^{k}$ be the $n\times n$ diagonal matrix obtained from the identity matrix $I_n$ by replacing the $(k,k)$-entry with $-1$. 
For a matrix $B=(b_{ij})$, let $[B]_+$ be the matrix obtained from $B$ by replacing every entry $b_{ij}$ with $[b_{ij}]_+$. Let $B^{k\bullet}$ ( resp. 
 $B^{\bullet k}$ ) be the matrix obtained from $B$ by replacing all entries outside of the $k$-th row(resp. column) with zeros.

\begin{proposition}
    The $F$-pattern associated with $(\mathbf{r},\mathbf{z})$-cluster pattern $\mathbf{\Sigma}$ is determined by the following initial condition and mutation formula:
    \begin{align}
        &F_{t_0}=\mathbf{0}_{n\times n},\\
        \label{f_vector_recurrence}
    &F_{t'}^{B_{t_0};t_0}=F_{t}^{B_{t_0};t_0}(J_{n}^k+[-\varepsilon B_{t}R]_+^{\bullet k})+[-\varepsilon C_{t}^{B_{t_0};t_0}R]_{+}^{\bullet k}+[\varepsilon C_{t}^{-B_{t_0};t_0}R]_{+}^{\bullet k},
    \end{align}
where $t$ and $t'$ are $k$-adjacent, $B_t$ is the principal part of $\tilde{B}_t$.
\end{proposition}
\begin{proof}
    The recurrence formulas for $F$-polynomials can be rewritten as 
\begin{align*}
    F_{k;t'}(\mathbf{\hat{y}}_{t_0},\mathbf{z})&=  F_{k;t}(\mathbf{\hat{y}}_{t_0},\mathbf{z})^{-1}\sum\limits_{s=0}^{r_{k}}z_{k,s}(\prod\limits_{j=1}^{n}\hat{y}_{j;t_0}^{s\varepsilon c_{jk;t}^{B_{t_0};t_0}+r_{k}[-\varepsilon c_{jk;t}^{B_{t_0};t_0}]_{+}}F_{j;t}(\mathbf{\hat{y}}_{t_0},\mathbf{z})^{s\varepsilon b_{jk;t}+r_{k}[-\varepsilon b_{jk;t}]_{+}})\\
   &= F_{k;t}(\mathbf{\hat{y}}_{t_0},\mathbf{z})^{-1}\sum\limits_{s=0}^{r_{k}}z_{k,s}(\hat{\mathbf{y}}_{t_0}^{s[\varepsilon \mathbf{c}_{k;t}^{B_{t_0};t_0}]_+}(\prod\limits_{j=1}^{n}F_{j;t}(\mathbf{\hat{y}}_{t_0},\mathbf{z})^{[\varepsilon b_{jk;t}]_+})^s\\
   &\quad\quad\hat{\mathbf{y}}_{t_0}^{(r_{k}-s)[-\varepsilon \mathbf{c}_{k;t}^{B_{t_0};t_0}]_{+}}(\prod\limits_{j=1}^{n}F_{j;t}(\mathbf{\hat{y}}_{t_0},\mathbf{z})^{[-\varepsilon b_{jk;t}]_{+}})^{r_k-s}).
\end{align*}
Since $ F_{k;t'}(\mathbf{\hat{y}}_{t_0},\mathbf{z})$ and $ F_{k;t}(\mathbf{\hat{y}}_{t_0},\mathbf{z})$ have ``maximal degree" monomials, it follows that
\[\sum\limits_{s=0}^{r_{k}}z_{k,s}(\hat{\mathbf{y}}_{t_0}^{s[\varepsilon \mathbf{c}_{k;t}^{B_{t_0};t_0}]_+}(\prod\limits_{j=1}^{n}F_{j;t}(\mathbf{\hat{y}}_{t_0},\mathbf{z})^{[\varepsilon b_{jk;t}]_+})^s\hat{\mathbf{y}}_{t_0}^{(r_{k}-s)[-\varepsilon \mathbf{c}_{k;t}^{B_{t_0};t_0}]_{+}}(\prod\limits_{j=1}^{n}F_{j;t}(\mathbf{\hat{y}}_{t_0},\mathbf{z})^{[-\varepsilon b_{jk;t}]_{+}})^{r_k-s})\] has a unique monomial divided by the other monomials in it, whose exponent is
\begin{eqnarray*}
&&\max\limits_{0\leq s\leq r_k}\{s([\varepsilon\mathbf{c}_{k;t}^{B_{t_0};t_0}]_++\sum_{j=1}^n[\varepsilon b_{jk;t}]_+\mathbf{f}_{j;t}^{t_0})+(r_k-s)([-\varepsilon\mathbf{c}_{k;t}^{B_{t_0};t_0}]_++\sum_{j=1}^n[-\varepsilon b_{jk;t}]_+\mathbf{f}_{j;t}^{t_0})\}\\
&=&\max\{r_k([\varepsilon\mathbf{c}_{k;t}^{B_{t_0};t_0}]_++\sum_{j=1}^n[\varepsilon b_{jk;t}]_+\mathbf{f}_{j;t}^{t_0}),r_k([-\varepsilon\mathbf{c}_{k;t}^{B_{t_0};t_0}]_++\sum_{j=1}^n[-\varepsilon b_{jk;t}]_+\mathbf{f}_{j;t}^{t_0})\}.
\end{eqnarray*}
Consequently, the recurrence formula of $f$-vectors is
\begin{align*}
    \mathbf{f}_{k;t'}^{t_{0}}&=
        -\mathbf{f}_{k;t}^{t_{0}}+\operatorname{max}\{r_{k}([\varepsilon\mathbf{c}_{k;t}^{B_{t_0};t_0}]_{+}+\sum\limits_{j=1}^{n}[\varepsilon b_{jk;t}]_{+}\mathbf{f}_{j;t}^{t_{0}}),r_{k}([-\varepsilon\mathbf{c}_{k;t}^{B_{t_0};t_0}]_{+}+\sum\limits_{j=1}^{n}[-\varepsilon b_{jk;t}]_{+}\mathbf{f}_{j;t}^{t_{0}})\}\\
        &=
        -\mathbf{f}_{k;t}^{t_{0}}+r_{k}([-\varepsilon\mathbf{c}_{k;t}^{B_{t_0};t_0}]_{+}+\sum\limits_{j=1}^{n}[-\varepsilon b_{jk;t}]_{+}\mathbf{f}_{j;t}^{t_{0}})+\operatorname{max}\{r_{k}(\varepsilon\mathbf{c}_{k;t}^{B_{t_0};t_0}+\sum\limits_{j=1}^{n}\varepsilon b_{jk;t}\mathbf{f}_{j;t}^{t_{0}}),0\}\\
        &=
        -\mathbf{f}_{k;t}^{t_{0}}+r_{k}([-\varepsilon\mathbf{c}_{k;t}^{B_{t_0};t_0}]_{+}+\sum\limits_{j=1}^{n}[-\varepsilon b_{jk;t}]_{+}\mathbf{f}_{j;t}^{t_{0}})+r_{k}[(\varepsilon \mathbf{c}_{k;t}^{B_{t_0};t_0}+\sum\limits_{j=1}^{n}\varepsilon b_{jk;t}\mathbf{f}_{j;t}^{t_{0}})]_+\\
        &=  -\mathbf{f}_{k;t}^{t_{0}}+r_{k}([-\varepsilon\mathbf{c}_{k;t}^{B_{t_0};t_0}]_{+}+\sum\limits_{j=1}^{n}[-\varepsilon b_{jk;t}]_{+}\mathbf{f}_{j;t}^{t_{0}})+r_{k}[\varepsilon \mathbf{c}_{k;t}^{-B_{t_0};t_0}]_+,
\end{align*}
where the last equality follows \eqref{eqn:c-matrix-duality}.
\end{proof}

\subsection{Initial-seed mutation formula of $F$-matrices}
In this section, we first prove some fundamental properties of $f$-vectors and then establish the initial-seed mutation formula for $F$-matrices. Let us fix another $(\bar{\mathbf{r}},\bar{\mathbf{z}})$-cluster pattern $\bar{\mathbf{\Sigma}}:=\{\bar{\Sigma}_t=(\bar{\mathbf{x}}_t,\bar{B}_t)\}_{t\in \mathbb{T}_n}$. In particular, $\bar{B}_t\in M_{m'\times n}(\mathbb{Z})$ for some positive integer $m'\geq n$. For each $t\in \mathbb{T}_n$,
denote by $\bar{B}_t^{pr}$ the principal part of $\bar{B}_t$ and
$\bar{\mathbf{x}}_t=(\bar{x}_{1;t},\dots, \bar{x}_{l;t})$. For each $i\in [1,n]$, denote by $\bar{\mathbf{f}}_{i;t}^{t_0}$ the $i$th $f$-vector and $\bar{F}_t^{t_0}$ the $F$-matrix of $\bar{\mathbf{\Sigma}}$ at vertex $t\in \mathbb{T}_n$ with respect to $t_0$.

\begin{proposition}\label{relation_between_f-pattern_and_B-pattern}
 Assume that $B_{t_0}R=\bar{B}_{t_0}^{pr}\bar{R}$, then $R^{-1}\mathbf{f}_{i;t}^{t_0}=\bar{R}^{-1}\mathbf{\bar{f}}_{i;t}^{t_0}$ for any $t\in \mathbb{T}_n$ and $i\in [1,n]$.
 
\end{proposition}
\begin{proof}
Let $\bar{C}_t^{\bar{B}^{pr}_{t_0};t_0}$ be the $C$-matrix of $\bar{\mathbf{\Sigma}}$ at vertex $t$ and  $\bar{\mathbf{c}}_{i;t}$ its $i$th column vector. By Proposition \ref{relation_of_C-matriice}, we have 
$R^{-1}C_t^{B_{t_0};t_0}R=\bar{R}^{-1}\bar{C}_{t}^{\bar{B}_{t_0}^{pr},t_0}\bar{R}$ and $R^{-1}C_t^{-B_{t_0},t_0}R=\bar{R}^{-1}\bar{C}_t^{-\bar{B}_{t_0}^{pr};t_0}R$. 
By left-multiplying $R^{-1}$ on both side of equation (\ref{f_vector_recurrence}), we conclude that $R^{-1}F_t^{B_{t_0};t_0}$ only depends on $B_{t_0}R$, $t_0$ and $t$. It follows that $R^{-1}F_{t}^{B_{t_0};t}=\bar{R}^{-1}\bar{F}_t^{\bar{B}_{t_0}^{pr};t_0}$.

\end{proof}
\begin{remark}
According to the proof of Proposition \ref{relation_between_f-pattern_and_B-pattern}, $R^{-1}\mathbf{f}_{i;t}^{t_{0}}$ is the $i$th $f$-vector at $t$ for the $B$-pattern $\{B_{t}R\}_{t\in\mathbb{T}_{n}}$.
    Similarly, $r_i^{-1}\mathbf{f}_{i;t}^{t_0}$ is the $i$th $f$-vector at $t$ for the $B$-pattern $\{RB_t\}_{t\in\mathbb{T}_n}$.
\end{remark}
The following results are useful for studying properties of $f$-vectors for generalized cluster algebras.
\begin{theorem}\label{two generalized coincide in cluster}
\cite[Theorem 3.7]{Cao_Li21} If $B_{t_{0}}R=\bar{B}_{t_{0}}^{pr}\bar{R}$, then the following statements hold:
\begin{enumerate}
    \item For any two vertices $t,t'\in\mathbb{T}_n$ and $1\le i,j\le n$, $x_{i;t}=x_{j;t'}$ if and only if $\bar{x}_{i;t}=\bar{x}_{j;t'}$;
    \item There exists a bijection $\alpha:\mathcal{X}(\mathbf{\Sigma})\to\mathcal{X}(\bar{\mathbf{\Sigma}})$ given by $\alpha(x_{i;t})=\bar{x}_{i;t}$. 
\end{enumerate}
    
\end{theorem}
\begin{theorem}\label{f-vector and cluster variables}
    \cite[Theorem 3.3]{Fu_Gyoda} 
    Assume that $\bar{\mathbf{r}}=(1,\dots, 1)$ and $\bar{\mathbf{z}}=\varnothing$. In particular, $\bar{\mathbf{\Sigma}}$ is an ordinary cluster pattern.
    The following statements hold:
    \begin{enumerate}
        \item The equality $\bar{f}_{ij;t'}^{t}=\bar{f}_{kl;s'}^{s}$ holds if $\bar{x}_{i;t}=\bar{x}_{k;s}$ and $\bar{x}_{j;t'}=\bar{x}_{l;s'}$, where $\bar{f}_{ij;t'}^{t}$ is the $i$th entry of $\bar{\mathbf{f}}_{j;t'}^{t}$ and  $\bar{f}_{kl;s'}^{s}$ is the $k$th entry  of $\bar{\mathbf{f}}_{l;s'}^s$.
            \item There exists a cluster containing $\bar{x}_{i;t}$ and $\bar{x}_{k;t_{0}}$ if and only if $\bar{f}_{ki;t}^{t_{0}}=0$.
            \item The cluster variable $x_{i;t}$ belongs to $\mathbf{x}_{t_0}$ if and only if $\mathbf{f}_{i;t}^{t_0}=0$.
    \end{enumerate}
\end{theorem}
Now we can state  the basic properties of $f$-vectors for generalized cluster algebras. As in \cite{Fu_Gyoda}, we refer to Theorem \ref{relation of f-compatible and cluster} $(2)$ as the {\it symmetry property}, and $(3)$ as the {\it compatibility property}. The proof of Theorem \ref{relation of f-compatible and cluster} (2) will be postponed  until after Proposition \ref{the relation of F-invariant and f-vector}.

\begin{theorem}
\label{relation of f-compatible and cluster}
Let $\mathbf{F}=\{F_t^{t_0}=(\mathbf{f}_{1;t}^{t_0},\dots, \mathbf{f}_{n;t}^{t_0})\}_{t\in \mathbb{T}_n}$ be the $F$-pattern associated with the $(\mathbf{r},\mathbf{z})$-cluster pattern $\mathbf{\Sigma}$. The following statements hold:
        \begin{enumerate}
            \item The equality $f^{t}_{ij;t'}r_{i}^{-1}=f^{s}_{kl;s'}r_{k}^{-1}$ holds if $x_{i;t}=x_{k;s}$ and $x_{j;t'}=x_{l;s'}$.
            \item  $DF_{t'}^tD^{-1}=(F_{t}^{t'})^\top$ for any $t,t'\in \mathbb{T}_n$.
            \item There exists an $(\mathbf{r},\mathbf{x})$-cluster containing $x_{i;t}$ and $x_{k;t_{0}}$ if and only if $f^{t_0}_{ki;t}=0$.
            \item The $(\mathbf{r},\mathbf{z})$-cluster variable $x_{i;t}$ belongs to $\mathbf{x}_{t_0}$ if and only if $\mathbf{f}_{i;t}^{t_0}=0$.\label{criteria of initial cluster using f-vector}

        \end{enumerate}
    \end{theorem}
    \begin{proof}
    Assume that $\bar{\mathbf{\Sigma}}=\{\bar{\Sigma}_t=(\bar{\mathbf{x}}_t,\bar{B}_t)\}_{t\in \mathbb{T}_n}$ is an ordinary cluster pattern, that is, $\bar{\mathbf{r}}=(1,\dots, 1)$ and $\bar{\mathbf{z}}=\varnothing$. Assume that $B_{t_{0}}R=\bar{B}_{t_{0}}^{pr}\bar{R}$.
       
           Suppose that $x_{i;t}=x_{k;s}$ and $x_{j;t'}=x_{l;s'}$ for $i,j,k,l\in[1,n]$, vertices $t,t',s,s'\in\mathbb{T}_n$.
 By Theorem \ref{two generalized coincide in cluster}, we have $\bar{x}_{i;t}=\bar{x}_{k;s}$ and $\bar{x}_{j;t'}=\bar{x}_{l;s'}$. It follows that $\bar{f}_{ij;t'}^{t}=\bar{f}^{s}_{kl;s'}$  by Theorem \ref{f-vector and cluster variables}. Thus, equality $f^{t}_{ij;t'}r_{i}^{-1}=f^{s}_{kl;s'}r_{k}^{-1}$ holds by Proposition \ref{relation_between_f-pattern_and_B-pattern}. This completes the proof of $(1)$. 
 
 The statements $(3)$ and $(4)$ follow directly from from Proposition \ref{relation_between_f-pattern_and_B-pattern}, Theorem \ref{two generalized coincide in cluster} and Theorem \ref{f-vector and cluster variables}.

We will prove $(2)$ using the $F$-invariant in the next section.
 \end{proof}


Recall that $\{C_t\mid t\in \mathbb{T}_n\}$ be the $(\mathbf{r},\mathbf{z})$-$C$-pattern associated with the $(\mathbf{r},\mathbf{z})$-cluster pattern $\mathbf{\Sigma}=\{\Sigma_t=(\mathbf{x}_t,\tilde{B}_t)\}_{t\in\mathbb{T}_n}$ with respect to $t_0\in \mathbb{T}_n$. The following is a direct consequence of the recurrence \eqref{E9}.
\begin{lemma}\label{duality of C-matrix}
    $\{DC_{t}D^{-1}|t\in \mathbb{T}_n\}$ is the ordinary $C$-pattern associated with the ordinary  $B$-pattern $\{-RB^{\top}_t|t\in \mathbb{T}_n\}$.
\end{lemma}

By Theorem \ref{relation of f-compatible and cluster} and the tropical duality between $C$-matrices and $G$-matrices in \cite{Nakanishi11}, we can obtain the initial-seed mutations of $F$-matrices for generalized cluster algebras.
\begin{proposition}\label{prop:initial-seed-mutation}
    Let $\{G_{s}^{B_tR;t}\}_{s\in \mathbb{T}_n}$ and $\{G_{s}^{-B_tR;t}\}_{t\in \mathbb{T}_n}$ be the ordinary $G$-pattern respectively associated with the ordinary $B$-pattern $\{B_sR\}_{s\in \mathbb{T}_n}$ and $\{-B_sR\}_{s\in \mathbb{T}_n}$ with respect to $t\in \mathbb{T}_n$. For $\xymatrix{t\ar@{-}[r]^k&t'}\in \mathbb{T}_n$, the following statement holds:
 \begin{align}
     F_{t_0}^{t'}=(J_n^k+[-\varepsilon RB_t]^{k\bullet}_+)F_{t_0}^{t}+r_k[\varepsilon G_{t_0}^{-B_{t}R;t}]_+^{k\bullet}+r_k[-\varepsilon G_{t_0}^{B_{t}R;t}]_+^{k\bullet}.\label{final seed mutation of f vector}
 \end{align}
\end{proposition}
\begin{proof}
According to Theorem \ref{relation of f-compatible and cluster} (2),  we have $DF_{s}^{s'}D^{-1}=(F_{s'}^s)^\top$ for any $s,s'\in \mathbb{T}_n$.
By \eqref{f_vector_recurrence}, we obtain
\begin{eqnarray*}
    F_{t_0}^{t'}&=&D^{-1}(F_{t'}^{t_0})^\top D\\
    &=&D^{-1}(F_{t}^{t_0}(J_n^k+[-\varepsilon B_tR]_+^{\bullet k}))+[-\varepsilon C_t^{B_{t_0};t_0}R]_+^{\bullet k}+[\varepsilon C_t^{-B_{t_0};t_0}R]_+^{\bullet k})^{\top}D\\
    &=&D^{-1}(J_n^k+[-\varepsilon B_tR]_+^{\bullet k})^{\top}D D^{-1}(F_t^{t_0})^{\top}D+D^{-1}([-\varepsilon C_t^{B_{t_0};t_0}R]_+^{\bullet k})^\top [-(\varepsilon C_t^{-B_{t_0};t_0}R]_+^{\bullet k})^\top D\\
    &=&(J_n^k+[\varepsilon RB_t]_{+}^{k\bullet})F_{t_0}^t+([-\varepsilon DC_t^{B_{t_0};t_0}D^{-1}R]_+^{\bullet k})^\top +([\varepsilon DC_t^{-B_{t_0};t_0}D^{-1}R]_+^{\bullet k})^\top\\
    &=&(J_n^k+[\varepsilon RB_t]_{+}^{k\bullet})F_{t_0}^t+r_k[-\varepsilon G_{t_0}^{-B_tR;t}]_+^{k\bullet}+r_k[\varepsilon G_{t_0}^{B_tR;t}]_+^{k\bullet},
\end{eqnarray*}
where the last equality follows from Lemma \ref{duality of C-matrix} and \cite[(1.13)]{Nakanishi11}.

\end{proof}
\section{F-invariant in generalized cluster algebras}\label{F-invariant in generalized cluster algebras}
\subsection{Tropical points and good elements}\label{tropical points and good elements}
In this section, we 
extend the definition of $F$-invariant \cite{cao23} to generalized cluster algebras and prove the symmetry property of $F$-matrix (Proposition \ref{relation of f-compatible and cluster} (2)).
Fix a Langland-Poisson triple $(\mathcal{S}_{X},\mathcal{S}_{Y},\Lambda)$ throughout this section. 
Let $\mathcal{A}$ be the $(\mathbf{r},\mathbf{z})$-cluster algebra associated with $\mathcal{S}_X$ and $\mathcal{U}$ be the upper $(\mathbf{r},\mathbf{z})$-cluster algebra associated with $\mathcal{S}_X$. 
Fix two {\it tropical semi-fields} $\mathbb{Z}^{\min}:=(\mathbb{Z},+,\oplus=\operatorname{min}\{-,-\})$ and $\mathbb{Z}^{\max}=(\mathbb{Z},+,\oplus=\operatorname{max}\{-,-\})$. For any semi-field $\mathbb{P}$ and $\mathbf{p}=(p_1,\dots,p_m)$, there exists a unique semi-field homomorphism
\begin{align*}
    \pi_{\mathbf{P}}:\mathcal{F}_{>0}&\to \mathbb{P}.\\
    z_{i,s}&\mapsto \id\\
    x_{i}&\mapsto p_{i}
\end{align*}

The map $\mathbf{p}\mapsto \pi_{\mathbf{p}}$ induces a bijection from $\mathbb{P}^{m}$ to set $\text{Hom}_{ssf}(\mathcal{F}_{>0},\mathbb{P})$ of special semi-field homomorphisms from $\mathcal{F}_{>0}$ to $\mathbb{P}$ mapping $z_{i,s}$ to identity in $\mathbb{P}$, where $1\le i\le n,1\le s \le r_i-1$. The elements in $\text{Hom}_{ssf}(\mathcal{F}_{>0},\mathbb{Z}^{\max})$ are called {\it tropical $\mathbb{Z}^{\max}$-points} or simply {\it tropical points} (cf. \cite[Section 3.1]{cao23}). 
A {\it chart} is an  $m$-tuple $\mathbf{u}=(u_1,\dots, u_n)$ of elements in $\mathcal{F}_{>0}$ such that $(\mathbf{z},\mathbf{u})$ generated $\mathcal{F}_{>0}$ as a semifield.
It is known that $\mathbf{x}_{t}=(x_{1;t},\cdots,x_{m;t})$ is a chart on $\mathcal{F}_{>0}$ for each vertex $t\in \mathbb{T}_n$. The following is proved in \cite{cao23}.
\begin{proposition}\cite[Proposition 3.2]{cao23}
    Let $\mathcal{C}:=\{\mathbf{u}_{t}|t\in\mathbb{T}_n\}$ be a collection of charts on $\mathcal{F}_{>0}$ indexed by set $\mathbb{T}_n$. Then the following statements hold.
    \begin{enumerate}
        \item For any chart $\mathbf{u}_{t}=(u_{1;t},\dots,u_{m;t})$, $x_{i}(\mathbf{u}_{t})$ denote the expression of $x_i$ in $\mathcal{F}_{>0;t}:=\mathbb{Q}_{sf}(\mathbf{z},\mathbf{u}_{t})$. There exists a map $\phi_{t}:x_i\mapsto x_{i}(\mathbf{u}_{t})$ inducing an isomorphism form $\mathcal{F}_{>0}$ to $\mathcal{F}_{>0;t}:=\mathbb{Q}_{sf}(\mathbf{z},\mathbf{u}_{t})$.
        \item Let $\nu$ be a tropical point in $\Hom_{ssf}(\mathcal{F}_{>0},\mathbb{Z}^{\max})$ and $\mathbf{u}_{t}=(u_{1;t},\dots,u_{m;t})$ be a chart on $\mathcal{F}_{>0}$. Set $q_{t}(\nu):=(\nu(u_{1;t}),\dots,\nu(u_{m;t}))\in\mathbb{Z}^{m}$. There exists a unique semi-field homomorphism $\pi_{q_{t}(\nu);t}$ in $\Hom_{ssf}(\mathcal{F}_{>0},\mathbb{Z}^{\max})$ induced by $u_{i;t}\mapsto \nu(u_{i;t})$ for any $i$ such that the following diagram commutes
$$\xymatrix{\mathcal{F}_{>0}\ar[r]^{\phi_{t}} \ar[d]^{\nu}&\mathcal{F}_{>0;t}\ar[dl]^{\pi_{q_{t}(\nu);t}}\\
\mathbb{Z}^{\max}&}$$
\item   The map $\nu\mapsto q_{t}(\nu)$ gives a bijection between $\text{Hom}_{ssf}(\mathcal{F}_{>0},\mathbb{Z}^{\max})$ and $\mathbb{Z}^{m}$. 

    \end{enumerate}
\end{proposition}

The integer vector $q_{t}(\nu)$ is called the coordinate vector of $\nu$ under $\mathbf{u}_{t}$. When the collection $\mathcal{C}=\{\mathbf{u}_{t}|t\in\mathbb{T}_n\}$ of charts is given, we identify the tropical point $\nu\in\text{Hom}_{ssf}(\mathcal{F}_{>0},\mathbb{Z}^{\max})$ with the collection $\{q_{t}(\nu)|t\in\mathbb{T}_n\}$ of coordinate vectors of $\nu$.
If we concentrate on the $Y$-pattern $\mathcal{S}_{Y}$, 
 the corresponding tropical points can be defined using the tropical version of the transition maps of $Y$-seed. 
\begin{proposition}
    \begin{enumerate}
        \item  A collection $[\mathbf{g}]:=\{\mathbf{g}_{t}\in\mathbb{Z}^{m}|t\in\mathbb{T}_{n}\}$ is a tropical point associated with $\mathcal{S}_{Y}$ if and only if 
 it satisfies the following recurrence relation:
\begin{align}
    g_{i;t'}=\begin{cases}
        -g_{i;t} &\quad if\ \ i=k.\\
        g_{i;t}+[r_{k}\hat{b}_{ki;t}]_{+}g_{k;t}+(-r_{k}\hat{b}_{ki;t})[g_{k;t}]_{+}&\quad if \ \ i\neq k,
    \end{cases}\label{recurrence of tropical points in Y}
\end{align}
for each edge $\xymatrix{t\ar@{-}[r]^{k}&t'}$ in $\mathbb{T}_{n}$. We denote by $\mathcal{S}_{Y}(\mathbb{Z}^{\max})$ the set of all tropical points associated with $\mathcal{S}_{Y}$.
\item  A collection $[\mathbf{a}]:=\{\mathbf{a}_{t}\in\mathbb{Z}^{m}|t\in\mathbb{T}_{n}\}$ is a tropical point associated with $\mathcal{S}_{X}$ if and only if 
 it satisfies the following recurrence relation:
\begin{align}
    a_{i;t'}=\begin{cases}
       a_{i;t} &\quad if \ \ i\neq k,\\
       -a_{k;t}+\max\{\sum\limits^{m}_{j=1}[b_{jk;t}r_{k}]_{+}a_{j;t},\sum\limits_{j=1}^{m}[-b_{jk;t}r_{k}]_{+}a_{j;t}\}&\quad if\ \ i=k,
    \end{cases}\label{recurrence of tropical points in X}
\end{align}
for each edge $\xymatrix{t\ar@{-}[r]^{k}&t'}$ in $\mathbb{T}_{n}$. We denote by $\mathcal{S}_{X}(\mathbb{Z}^{\max})$ the set of all tropical points associated with $\mathcal{S}_{X}$.
    \end{enumerate}    
\end{proposition}
 Let $\tilde{B}_{t_0}^{sq}$ be the skew-symmetrizable $m\times m$ matrix $(\tilde{B}_{t_0}\mid M)$. Denote $\{\tilde{B}_t^{sq}|t\in\mathbb{T}_n\}$ the collection of matrices obtained from $\tilde{B}_{t_0}^{sq}$ by any sequence of mutations in directions $1,\cdots,n$. Let  $\tilde{D} = \diag(d_1, \ldots, d_m)$ be the skew-symmetrizer of $\tilde{B}_{t}^{sq}$. There exist canonical duality maps between  $\mathcal{S}_{X}(\mathbb{Z}^{T})$ and $\mathcal{S}_{Y}(\mathbb{Z}^{T})$, which generalize the tropical duality in \cite[Proposition 4.3 and 4.4]{cao23} to generalized cluster algebras.
\begin{proposition}[Tropical Duality Correspondence]\label{maps between two tropical points}

    (1) There is a map $\Phi:\mathcal{S}_{X}(\Z^{\max})\to\mathcal{S}_{Y}(\Z^{\max})$ such that for any tropical point $[\mathbf{a}] = \{\mathbf{a}_t \in \Z^m \mid t \in \mathbb{T}_n\}$ in $\mathcal{S}_X(\Z^{\max})$, 
    \[
        \Phi([\mathbf{a}]) := \left\{\tilde{D} (\widetilde{B}_t^{\mathrm{sq}})^\top \mathbf{a}_t \in \Z^m \,\big|\, t \in \mathbb{T}_n \right\}.
    \]
    
   (2) There is a map $\Psi:\mathcal{S}_{Y}(\Z^{\max})\to\mathcal{S}_{X}(\Z^{\max})$ such that  for any tropical point $[\mathbf{g}] = \{\mathbf{g}_t \in \Z^m \mid t \in \mathbb{T}_n\}$ in $\mathcal{S}_Y(\Z^T)$, 
    \[
        \Psi([\mathbf{g}]) := \left\{ \Lambda_t \mathbf{g}_t \in \Z^m \,\big|\, t \in \mathbb{T}_n \right\}.
    \]

\end{proposition}
\begin{proof}
    The same proof of \cite[Proposition 4.3 and 4.4]{cao23} is valid.
\end{proof}

Assume that $\mathcal{U}$ is a full rank upper $(\mathbf{r},\mathbf{z})$-cluster algebra, that is, $\text{rank}(\tilde{B}_t)=n$. For each seed $(\mathbf{x}_{t},\tilde{B}_{t})$ of $\mathcal{U}$, the dominance partial order $\le_{t}$   on $\mathbb{Z}^{m}$ \cite{Qin17} is defined as follows:\\
For two vectors $\mathbf{g},\mathbf{g}'\in\mathbb{Z}^{m}$, we say $\mathbf{g}'\le_{t}\mathbf{g}$ if there exists a vector $\nu=(\nu_{1},\cdots,\nu_{n})^{\top}\in\mathbb{N}^{n}$ such that 
\[\mathbf{g}'=\mathbf{g}+\tilde{B}_{t}\nu.\]
We write $\mathbf{g}'\prec_{t}\mathbf{g}$ if $\mathbf{g}'\le_{t}\mathbf{g}$ and $\mathbf{g}'\neq\mathbf{g}$. 
For each $t\in \mathbb{T}_n$, denote by 
\[R_t:=\{\sum_{\mathbf{h}\in\mathbb{Z}^m}b_{\mathbf{h}}\mathbf{x}_{t}^{\mathbf{h}}\mid b_{\mathbf{h}}\in\mathbb{Z}[\mathbf{z}],\max_{\le_t}\{\mathbf{h}|b_{\mathbf{h}}\neq0\} \text{ is finite}\}\]
the ring of formal power Laurent series, 
where $\max_{\le_t}$ denotes the set of maximal elements under the dominance order. 
An element $u=\sum_{\mathbf{h}\in\mathbb{Z}^{m}}b_{\mathbf{h}}\mathbf{x}_{t}^{\mathbf{h}}\in R_{t}$ is called {\it pointed} for the seed $(\mathbf{x}_{t},\tilde{B}_{t})$ if :
\begin{enumerate}
    \item The support $\supp(u):=\{\mathbf{h}\in\mathbb{Z}^{m}|b_{\mathbf{h}}\neq0\}$ contains a unique maximal element $\mathbf{g}$ under $\le_{t}$, called the dominance degree $\deg^{t}(u)=\mathbf{g}$.
    \item The leading coefficient satisfies  $b_{\mathbf{g}}=1$.
\end{enumerate}
Every pointed element $u$ relative to $(\mathbf{x}_{t},\tilde{B}_{t})$ of $\mathcal{U}$ admits a canonical decomposition
\[u=\mathbf{x}_{t}^{\mathbf{g}}+\sum\limits_{\mathbf{h}\le_{t}\mathbf{g}}b_{\mathbf{h}}\mathbf{x}_{t}^{\mathbf{h}}=\mathbf{x}_{t}^{\mathbf{g}}F(\hat{y}_{1;t},\cdots,\hat{y}_{n;t}),\]
where $b_{\mathbf{h}}\in\mathbb{Z}[\mathbf{z}]$ and $F\in\mathbb{Z}[\mathbf{z}][\hat{y}_{1;t},\cdots,\hat{y}_{n;t}]$ is a polynomial with constant term $1$. Notably, all cluster monomials are pointed elements with respect to $(\mathbf{x}_t,\tilde{B}_t)$. We refer to \cite{cao23} for structural properties of $R_t$.
Analogous to classical cluster algebras, these pointed elements are governed by tropical points associated with the $Y$-pattern.
\begin{definition}\label{definition of good element}
    Let $(\mathcal{S}_{\mathbf{X}},\mathcal{S}_{\mathbf{Y}})$ be Langlands dual pair,  let $\mathcal{U}$ be full rank upper $(\mathbf{r},\mathbf{z})$-cluster algebra associated with $\mathcal{S}_{\mathbf{X}}$.
 \begin{enumerate}
       \item An element $u\in\mathcal{U}$ is compatibly pointed if:
\begin{itemize}
\item u is pointed for every seed $(\mathbf{x}_t, \tilde{B}_t)$;
\item   The collection $[\mathbf{g}]:=\{\deg^t(u)|t\in\mathbb{T}_n\}$ forms a tropical point in $\mathcal{S}_{\mathbf{Y}}(\mathbb{Z}^{\max})$. 
\end{itemize}
Such elements are termed $[\mathbf{g}]$-pointed, consistent with the framework in \cite{cao23}.
\item A $[\mathbf{g}] $-pointed element $ u \in \mathcal{U} $ is called $ [\mathbf{g}] $-good if it admits a subtraction-free rational expression in $ \mathbf{x}_t$ for every vertex $t \in \mathbb{T}_n$.
\end{enumerate}
\end{definition}
\begin{remark}
   Unlike \cite{cao23}, we omit the universal positivity requirement for good elements.
\end{remark}
Let $u=\mathbf{x}^{\mathbf{g}_{t}}_{t}F(\hat{y}_{1;t},\cdots,\hat{y}_{n;t})$ be $[\mathbf{g}]$-pointed element with $\mathbf{g}_{t}=\deg^{t}(u)$ and $F(\hat{y}_{1;t},\cdots,\hat{y}_{n;t})\in\mathbb{Z}[\mathbf{z}][\hat{y}_{1;t},\cdots,\hat{y}_{n;t}]$ having constant term $1$. 
We have the following notion:
\begin{enumerate}
    \item \textbf{$F$-Polynomial}: The polynomial $F(\hat{y}_{1;t},\cdots,\hat{y}_{n;t})$ is called the {\it $F$-polynomial} of $u$ at vertex \textit{t}.
    \item \textbf{$f$-Vector}: For $k\in[1,n]$, let $f_{k}^t(u):=\operatorname{max}\{d\in\mathbb{N}\mid \hat{y}_{k;t}^d\ \text{divides some terms in } F\}$. The $f$-vector of $u$ at $t$ is:\[\mathbf{f}^{t}=(f_{1}^t(u),\cdots,f_{n}^t(u))^{\top}\in\mathbb{N}^{n}\]
    \item \textbf{Bipointed Element}: $u$ is called $[\mathbf{g}]$-bipointed if for every $t\in\mathbb{T}_n$, the monomial $\prod_{j=1}^{n}\hat{y}_{j;t}^{f_{j}^t(u)}$ appears in $F(\hat{y}_{1;t},\cdots,\hat{y}_{n;t})$ with coefficient 1.
    \item \textbf{Bigood Element}: A $[\mathbf{g}]$-good element is $[\mathbf{g}]$-bigood if it is $[\mathbf{g}]$-bipointed.
\end{enumerate}
 
\begin{proposition}{\rm (\cf \cite[Proposition 3.7]{cao23})}
    Every $(\mathbf{r},\mathbf{z})$-cluster variable $x_{i;t}$ is $[\mathbf{g}]$-bigood, where $[\mathbf{g}]$ corresponds to the tropical point in $\mathcal{S}_Y$ defined by the extended $g$-vectors:
    $$[\mathbf{g}]:=\{\mathbf{g}_{i;t}^{w}\mid\forall w\in\mathbb{T}_n,\text{ $g_{i;t}^w$ is extended $g$-vector of $x_{i;t}$ w.r.t $\mathbf{x}_w$}\}.$$
    Consequently, all cluster monomials are bigood.
\end{proposition}
\begin{proof}
Let $\{G_t^{\tilde{B}_{t_0};t_0}\mid t\in \mathbb{T}_n\}$ be the extended $G$-pattern associated with the $(\mathbf{r},\mathbf{z})$-cluster pattern $\mathcal{S}_X$ with respect to $t_0$.
By the recurrence relation shown in \cite[(2.14)]{Fu_Peng_Ye}, for each edge $\xymatrix{t\ar@{-}[r]^{k}&t'}$ in $\mathbb{T}_{n}$, we have \[\tilde{G}^{\tilde{B}_{t_0};t_0}_{t'}=\tilde{G}_{t}^{\tilde{B}_{t_0};t_0}(J_m^k+[-\varepsilon_{k;t}\tilde{B}_tR]_+^{\bullet k}),\]
  where $\varepsilon_{k;t}$ is the common sign of entries of the $k$th row of $\tilde{G}_{t}^{\tilde{B}_{t_0};t_0}$. Let $\bar{B}_{t_0}=[\tilde{B}_{t_0}R\mid Q]$ be a skew-symmetry $m\times m$ matrix. Let $\{\bar{B}_{t}\mid t\in\mathbb{T}_m\}$ be the ordinary $B$-pattern associated with $\bar{B}_{t_0}$. Let $\{{G}^{\bar{B}_{t_0};t_0}_t\mid t\in \mathbb{T}_m\}$ be the ordinary $G$-pattern associated with ordinary $B$-pattern $\{\bar{B}_{t}\mid t\in\mathbb{T}_m\}$. By viewing $\mathbb{T}_n$ as a subgraph of $\mathbb{T}_m$, 
  we obtain $\tilde{G}^{\tilde{B}_{t_0};t_0}_{t}=G_t^{\bar{B}_{t_0};t_0}$, where $t$ and $t_0$ are linked by a sequence of edges labeled by $i\in[1,n]$. 
  By \cite[Proposition 3.6]{Fujiwara_Gyoda19}, we conclude that
  \[\tilde{G}^{\tilde{B}_{w'};w'}_{t}=(J_m^k+[\varepsilon_k\tilde{B}_wR]_+^{\bullet k})\tilde{G}^{\tilde{B}_w;w}_t,\]
  where $t$ and $\xymatrix{w\ar@{-}[r]^{k}&w'}\in \mathbb{T}_n$, and
  $\varepsilon_k$ is the common sign of entries of the $k$th row of $\tilde{G}_{t}^{\tilde{B}_{w};w}$.
  The result follows immediately.
\end{proof}
Applying almost the same proof of \cite[Corollary 3.14]{cao23} yields the following.

\begin{proposition}{\rm(cf. \cite[Corollary 3.14]{cao23})}\label{properties of f-vectors}
Let $u$ be a $[\mathbf{g}]$-pointed element in $\mathcal{U}$. Let $\mathbf{f}^{t}=(f_{1}^t,\cdots,f_{n}^t)^{\top}$ be the $f$-vector of $u$ with respect to $t\in\mathbb{T}_{n}$. Then the following statements hold:
    \begin{enumerate}
        \item If $f_{k}^t=0$ for some $k\in [1,n]$, then the kth component $g_{k;t}$ of $\mathbf{g}_{t}=\deg^{t}(u)$ is non-negative;
        \item If $f_{i}^t=0$ for any $i\in [1,n]$, then $u$ is a cluster monomial in $\mathbf{x}_{t}$
        \item Let $k\in [1,n]$ and $(\mathbf{x}_{t'},\tilde{B}_{t'})=\mu_{k}(\mathbf{x}_{t},\tilde{B}_{t})$. Suppose that $f_{k}^t\neq 0$ and $f_{i}^t=0$ for any $i\neq k$. Then $u$ is a cluster monomial in $\mathbf{x}_{t'}$.
    \end{enumerate}
\end{proposition}
\subsection{F-invariant}
%

We generalize the $F$-invariant \cite{cao23} to the setting of generalized cluster algebras, and prove that it has the same properties as it does in cluster algebras.

Recall that $(\mathcal{S}_X,\mathcal{S}_Y,\Lambda)$ is a {Langland-Poisson triple}, and  $\mathcal{U}$ is the upper $(\mathbf{r},\mathbf{z})$-cluster algebra associated with $\mathcal{S}_X$. By definition, we have
\[
\widetilde{B}_{t_0}^\top \Lambda_{t_0} = (D\mid 0),
\]
where $D = \operatorname{diag}(d_1, \ldots, d_n)$ and $0$ is the $n \times m$ zero matrix.

\begin{definition}
Let $\mathcal{U}^{[\mathbf{g}]}$ denote the set of $[\mathbf{g}]$-good elements in $\mathcal{U}$ for tropical points $[\mathbf{g}] \in \mathcal{S}_Y(\mathbb{Z}^{\max})$. The collection of all good elements is 
\[
\mathcal{U}^{\text{good}} := \bigcup_{[\mathbf{g}] \in \mathcal{S}_Y(\mathbb{Z}^T)} \mathcal{U}^{[\mathbf{g}]}.
\]
\end{definition}

\begin{notation}
For any non-zero rational function 
\[
F = \frac{\sum_{\nu \in \mathbb{N}^n} c_\nu \widehat{\mathbf{y}}^\nu}{\sum_{\mu \in \mathbb{N}^n} d_\mu \widehat{\mathbf{y}}^\mu},
\]
where $c_\nu, d_\mu \in \mathbb{Z}_{\geq 0}[\mathbf{z}]$ with $c_0 = d_0 = 1$ and, define for $\mathbf{h} \in \mathbb{Z}^n$:
\[
F[\mathbf{h}] := \max\{\nu^\top \mathbf{h} \mid c_\nu \neq 0\} - \max\{\mu^\top \mathbf{h} \mid d_\mu \neq 0\} \in \mathbb{Z}.
\]
\end{notation} 
    It is to see that $F[\mathbf{h}]$ is independent of the choice of the expression of $F$ as rational function by polynomials with non-negative coefficients. The following is obvious.
\begin{lemma}\label{lem:F-h}
    Let $F=\sum_{v\in \mathbb{N}^n}a_v\hat{\mathbf{y}}^v\in \mathbb{Z}[\mathbf{z}][\hat{\mathbf{y}}]$, $G=\sum_{v\in \mathbb{N}^n}c_v\hat{\mathbf{y}}^v, H=\sum_{v\in \mathbb{N}^n}d_v\hat{\mathbf{y}}^v\in \mathbb{Z}_{\geq 0}[\mathbf{z}][\hat{\mathbf{y}}]$  such that $F=G/H$ and $a_0=c_0=d_0=1$. Assume that $F$ has a unique monomial $\hat{\mathbf{y}}^f$ such that each monomial $\hat{\mathbf{y}}^v$ with nonzero coefficient in $F$ divides $\hat{\mathbf{y}}^f$ and the coefficient of $\hat{\mathbf{y}}^f$  is $1$. Then for any $\mathbf{h}\in \mathbb{N}^n$, $F[\mathbf{h}]=f^\top \mathbf{h}$.
\end{lemma}
\begin{definition}
    {\rm (F-invariant, cf. \cite[Section 4.2]{cao23})}\begin{enumerate}
        \item For any vertices $t\in\mathbb{T}_{n}$, we define a pair 
\begin{align*}
    \langle-\parallel-\rangle_{t}:\mathcal{U}^{good}\times \mathcal{U}^{good}\longrightarrow\mathbb{Z}
\end{align*}
by $\langle\mathbf{u}_{[\mathbf{g}]}\parallel\mathbf{u}_{[\mathbf{g'}]}\rangle_t=\mathbf{g}_{t}^{\top}\Lambda_{t}\mathbf{g'}_{t}+F_{t}[[D\mid0]\mathbf{g'}_{t}],$ where $\mathbf{u}_{[\mathbf{g}]}=\mathbf{x}_{t}^{\mathbf{g}_{t}}F_{t}(\hat{y}_{1;t},\cdots,\hat{y}_{n;t})$ and $\mathbf{u}_{[\mathbf{g'}]}=\mathbf{x}_{t}^{\mathbf{g'}_{t}}F'_{t}(\hat{y}_{1;t},\cdots,\hat{y}_{n;t})$;
\item Symmetrized version:
\begin{align*}
    (\mathbf{u}_{[\mathbf{g}]}\parallel\mathbf{u}_{[\mathbf{g'}]})_{t}&=\langle\mathbf{u}_{[\mathbf{g}]}\parallel\mathbf{u}_{[\mathbf{g'}]}\rangle_{t}+\langle\mathbf{u}_{[\mathbf{g'}]}\parallel\mathbf{u}_{[\mathbf{g}]}\rangle\\
    &=F_{t}[[D\mid 0]\mathbf{g'}_{t}]+F_{t}'[[D\mid0]\mathbf{g}_{t}].
\end{align*}
We will prove that $\langle-\parallel-\rangle_{t}$ is independent of the choice of vertex $t$ and we call $(\mathbf{u}_{[\mathbf{g}]}||\mathbf{u}_{[\mathbf{g'}]})_{t}$ the $F$-invariant of $\mathbf{u}_{[\mathbf{g}]}$ and $\mathbf{u}_{[\mathbf{g'}]}$.
\item Two pointed elements $\mathbf{u}_{[\mathbf{g}]}$ and $\mathbf{u}_{[\mathbf{g'}]}$ are said to be $F$-compatible, if 
\[(\mathbf{u}_{[\mathbf{g}]}\parallel\mathbf{u}_{[\mathbf{g'}]})_{t}=0.\]
    \end{enumerate}
\end{definition}
In particular, we obtain the following result, \cf \cite[Theorem 4.8]{cao23}.
\begin{theorem}\label{thm:independent-t}
    The pair $\langle-\parallel-\rangle_{t}$ is independent of the choice of vertex $t\in\mathbb{T}_{n}$. \end{theorem}
    \begin{proof}
       Let $\mathbf{u}_{[\mathbf{g}]}=\mathbf{x}_{t}^{\mathbf{g}_{t}}F_{t}(\hat{y}_{1;t},\cdots,\hat{y}_{n;t})$ and $\mathbf{u}_{[\mathbf{g'}]}=\mathbf{x}_{t}^{\mathbf{g'}_{t}}F'_{t}(\hat{y}_{1;t},\cdots,\hat{y}_{n;t})$ be $[\mathbf{g}]$-and $[\mathbf{g'}]$-good elements in $\mathcal{U}$ respectively. By definition, we have \[F_t=\frac{\sum\limits_{\nu\in\mathbb{N}^{n}}c_{\nu}\hat{\mathbf{y}}^{\nu}_t}{\sum\limits_{u\in\mathbb{N}^{n}}d_{u}\hat{\mathbf{y}}^{u}_t}\ \text{and}\ F_t'=\frac{\sum\limits_{\nu\in\mathbb{N}^{n}}c'_{\nu}\hat{\mathbf{y}}^{\nu}_t}{\sum\limits_{u\in\mathbb{N}^{n}}d'_{u}\hat{\mathbf{y}}^{u}_t},\]
       where $c_\nu,d_u,c_\nu',d_u'\in \mathbb{Z}_{\geq 0}[\mathbf{z}]$ and $c_0=d_0=c_0'=d_0'=1$. According to Proposition \ref{maps between two tropical points}, $\{\Lambda_{t}\mathbf{g'}_{t}\in\mathbb{Z}^{m}|t\in\mathbb{T}_{n}\}$ is a tropical point in $\mathcal{S}_{X}(\mathbb{Z}^{\max})$. There exist a morphism $\mu\in\operatorname{Hom}_{ssf}(\mathcal{F}_{>0},\mathbb{Z}^{\max})$ satisfying $\mu(x_{1;t},\cdots,x_{m;t})=\Lambda_{t}\mathbf{g'}_{t}$. Then, we have $\mu(\mathbf{x}_{t}^{\mathbf{g}_{t}})=\mathbf{g}_{t}^{\top}\Lambda_{t}\mathbf{g'}_{t}$ and \begin{align*}
           \mu(\mathbf{u}_{[\mathbf{g}]})&=\mu(\mathbf{x}_{t}^{\mathbf{g}_{t}}F_{t}(\hat{y}_{1;t},\cdots,\hat{y}_{n;t}))\\
           &=\mu(\mathbf{x}_{t}^{\mathbf{g}_{t}})+\mu(F_{t}(\hat{y}_{1;t},\cdots,\hat{y}_{n;t}))\\
           &=\mathbf{g}_{t}^{\top}\Lambda_{t}\mathbf{g'}_{t}+\mu(\sum\limits_{\nu\in\mathbb{N}^{n}}c_{\nu}\hat{\mathbf{y}}^{\nu}_t/\sum\limits_{u\in\mathbb{N}^{n}}d_{u}\hat{\mathbf{y}}^{u}_t)\\
           &=\mathbf{g}_{t}^{\top}\Lambda_{t}\mathbf{g'}_{t}+\operatorname{max}\{\mu(\hat{\mathbf{y}}_{t}^{\nu})|c_{\nu}\neq 0\}-\operatorname{max}\{\mu(\hat{\mathbf{y}}_{t}^{u})|d_{\mu}\neq 0\}\\
           &=\mathbf{g}_{t}^{\top}\Lambda_{t}\mathbf{g'}_{t}+\operatorname{max}\{(\tilde{B}_{t}\nu)^{\top}\Lambda_{t}\mathbf{g'}_{t}|c_{\nu}\neq 0\}-\operatorname{max}\{(\tilde{B}_{t}u)^{\top}\Lambda_{t}\mathbf{g'}_{t}|d_{u}\neq 0\}\\
&=\mathbf{g}_{t}^{\top}\Lambda_{t}\mathbf{g'}_{t}+\operatorname{max}\{\nu^{\top}\tilde{B}_{t}^{\top}\Lambda_{t}\mathbf{g'}_{t}|c_{\nu}\neq 0\}-\operatorname{max}\{u^{\top}\tilde{B}_{t}^{\top}\Lambda_{t}\mathbf{g'}_{t}|d_{u}\neq 0\}\\
&=\mathbf{g}_{t}^{\top}\Lambda_{t}\mathbf{g'}_{t}+\operatorname{max}\{\nu^{\top}[D\mid0]\mathbf{g'}_{t}|c_{\nu}\neq 0\}-\operatorname{max}\{u^{\top}[D\mid0]\mathbf{g'}_{t}|d_{u}\neq 0\}\\
&=\mathbf{g}_{t}^{\top}\Lambda_{t}\mathbf{g'}_{t}+F_{t}[[D\mid0]\mathbf{g'}_{t}]\\
&=\langle\mathbf{u}_{[\mathbf{g}]}\parallel\mathbf{u}_{[\mathbf{g'}]}\rangle_t
       \end{align*}
    Thus $\langle-\parallel-\rangle_{t}$ is independent of the choice of vertex $t\in\mathbb{T}_{n}$.
    \end{proof}
   As a consequence of Theorem \ref{thm:independent-t}, we can get rid of vertex $t$ and denote
    \begin{align*}
        \langle\mathbf{u}_{[\mathbf{g}]}\parallel\mathbf{u}_{[\mathbf{g'}]}\rangle_F:&=\langle\mathbf{u}_{[\mathbf{g}]}\parallel\mathbf{u}_{[\mathbf{g'}]}\rangle_t,\\
        (\mathbf{u}_{[\mathbf{g}]}\parallel\mathbf{u}_{[\mathbf{g'}]})_F:&=(\mathbf{u}_{[\mathbf{g}]}\parallel\mathbf{u}_{[\mathbf{g'}]})_t\\
        &=F_{t}[[D|0]\mathbf{g'}_t]+F'_{t}[[D|0]\mathbf{g}_{t}].
    \end{align*}
     In particular, we have $\langle x_{i;t}\parallel x_{j;t}\rangle_F=\langle x_{i;t}\parallel x_{j;t}\rangle_t=e_i^{\top}\Lambda_t e_{j}=\lambda_{ij;t}$ for any two cluster variables $x_{i;t}$ and $x_{j;t}$ in the same cluster.
     \begin{proposition}
         \label{the relation of F-invariant and f-vector}
         Let $x=x_{i;t}$ be an $(\mathbf{r},\mathbf{z})$-cluster variable for some $i\in [1,m]$ and $t\in\mathbb{T}_n$. Let $u$ be a $[\mathbf{g}]$-good element of $\mathcal{U}$. Let the Laurent expression of $u$ with respect to $t$ is \[u=\mathbf{x}_t^{\mathbf{g}_t}F(\hat{y}_{1;t},\cdots,\hat{y}_{n;t})=\mathbf{x}_t^{\mathbf{g}_t}\sum\limits_{\nu\in\mathbb{N}^{n}}a_{\nu}\hat{\mathbf{y}}_t^{\nu}\]and $\mathbf{f}^{t}=(f_{1}^t,\cdots,f_{n}^t)$ be $f$-vector of $u$ with respect to $t$. Then
         \begin{align*}
             (x\parallel u)_F=\begin{cases}
                 d_if_{i}^t&\quad if \ \ i\in [1,n],\\
                 0 &\quad if \ \ i\in [n+1,m],
             \end{cases}
         \end{align*} 
         where $d_i$ is the $(i,i)$-entry of $D=\operatorname{diag}(d_1,\cdots,d_n)$.
     \end{proposition}    \begin{proof}
         By the recurrence relation of the $F$-polynomial, we have $$F(\hat{y}_{1;t},\cdots,\hat{y}_{n;t})=\frac{\sum\limits_{\nu\in\mathbb{N}^{n}}c_{\nu}\hat{\mathbf{y}}^{\nu}}{\sum\limits_{u\in\mathbb{N}^{n}}d_{u}\hat{\mathbf{y}}^{u}},$$where $\sum\limits_{\nu\in\mathbb{N}^{n}}c_{\nu}\hat{\mathbf{y}}^{\nu},\sum\limits_{\mu\in\mathbb{N}^{n}}d_{\mu}\hat{\mathbf{y}}^{\mu}\in\mathbb{Z}_{>0}[\mathbf{z}][\hat{y}_{1},\cdots,\hat{y}_n]$ with constant term $1$. Let $f^1$ be the maximal degree of $\hat{y}_{i;t}$ in $\sum\limits_{\nu\in\mathbb{N}^{n}}c_{\nu}\hat{\mathbf{y}}^{\nu}$ and $f^2$ be the maximal degree of $\hat{y}_{i;t}$ in $\sum\limits_{\nu\in\mathbb{N}^{n}}d_{\mu}\hat{\mathbf{y}}^{\mu}$. We find $f^{1}-f^{2}=f_{i}^{t}$.
         Since the $F$-invariant is independent of the choice of vertex $t$. Then
         \begin{align*}
             (x\parallel u)_F&=(x_{i;t}||u)_t\\
             &=\operatorname{max}\{\nu^{\top}[D\mid0]e_{i}|c_{\nu}\neq 0\}-\operatorname{max}\{\mu^{\top}[D\mid0]e_{i}|d_{\mu}\neq 0\}\\
              &=\begin{cases}
                 d_if^{1}-d_if^{2}&\quad if \ \ i\in [1,n],\\
                 0 &\quad if \ \ i\in [n+1,m].
             \end{cases}\\
             &=\begin{cases}
                 d_if_{i}^t&\quad if \ \ i\in [1,n],\\
                 0 &\quad if \ \ i\in [n+1,m].
             \end{cases}
         \end{align*}
     \end{proof}
\begin{proof}[Proof of Theorem \ref{relation of f-compatible and cluster} (2).] By the symmetry of $F$-invariant, we get the symmetry of $f$-vectors
\begin{align*}
    d_if_{ij;t'}^{t}=(x_{i;t}\parallel x_{j;t'})_F=(x_{j;t'}\parallel x_{i;t})_F= d_jf_{ji;t}^{t'}.
\end{align*} 
\end{proof}

\begin{corollary}
     \label{criterion for the same cluster of cluster variables and good elements}
      Let  $u$ be a good element in $\mathcal{U}$. Let $\mathbf{x}_{t}$ be an $(\mathbf{r},\mathbf{z})$-cluster of $\mathcal{A}$. Then the following statements hold:
    \begin{enumerate}
        \item If $(x_{k;t}\parallel u)_F=0$ for some $k\in [1,n]$, then the kth component $g_{k;t}$ of $\mathbf{g}_{t}=\deg^{t}(u)$ is non-negative;
        \item If $(x_{i;t}\parallel u)_F=0$ for any $i\in [1,n]$, then $u$ is a cluster monomial in $\mathbf{x}_t$;
        \item Let $k\in [1,n]$ and $(\mathbf{x}_{t'},\tilde{B}_{t'})=\mu_{k}(\mathbf{x}_{t},\tilde{B}_{t})$. Suppose that $(x_{k;t}\parallel u)_F\neq 0$ and $(x_{i;t}\parallel u)_F=0$ for any $i\neq k$. Then $u$ is a cluster monomial in $\mathbf{x}_{t'}$.
        \end{enumerate}
\end{corollary}
\begin{proof}
    The results follow from Proposition \ref{the relation of F-invariant and f-vector} and Corollary \ref{properties of f-vectors}.
\end{proof}
\begin{proposition}\label{liner F-invariant}
    Let $u$ be a $[\mathbf{g}]$-bigood element in $\mathcal{U}$ and let $\mathbf{x}_t^{\mathbf{h}}=\prod\limits_{j=1}^{m}x_{j;t}^{h_j}$ be a cluster monomial of $\mathcal{U}$, where $\mathbf{h}=(h_1,\cdots,h_m)\in\mathbb{N}^m$. Then
    \begin{align*}
        (u\parallel\mathbf{x}_{t}^{\mathbf{h}})_F=\sum\limits_{j=1}^{m}h_{j}(u\parallel x_{j;t})_{F}.
    \end{align*}
\end{proposition}
\begin{proof}
    Since $u$ is $[\mathbf{g}]$-bigood element in $\mathcal{U}$, we have 
    \begin{align*}
        u&=\mathbf{x}_{t}^{\mathbf{g}_{t}}F(\hat{y}_{1;t},\cdots,\hat{y}_{n;t})\\
        &=\mathbf{x}_{t}^{\mathbf{g}_t}(\sum\limits_{\nu\in\mathbb{N}^{n}}c_{\nu}\hat{\mathbf{y}}_{t}^{\nu}/\sum\limits_{u\in\mathbb{N}^{n}}d_{u}\hat{\mathbf{y}}_{t}^{u}),
    \end{align*}
    where $\sum\limits_{\nu\in\mathbb{N}^{n}}c_{\nu}\hat{\mathbf{y}}_{t}^{\nu},\sum\limits_{\mu\in\mathbb{N}^{n}}d_{\mu}\hat{\mathbf{y}}_{t}^{\mu}\in\mathbb{Z}_{\geq 0}[\mathbf{z}][\hat{y}_{1},\cdots,\hat{y}_n]$ with constant term $1$. Let $\mathbf{f}=(f_1,\dots,f_n)^{\top}$ be the $f$-vector of $u$ with respect to $t$. By Lemma \ref{lem:F-h}, we obtain
     \begin{align*}
        (u\parallel\mathbf{x}_{t}^{\mathbf{h}})_F&=(u\parallel\mathbf{x}_{t}^{\mathbf{h}})_t\\
        &=F[[D|0]\mathbf{h}]\\
        &=\mathbf{f}^{\top}[D|0]\mathbf{h}\\
        &=\sum\limits_{i=1}^{n}d_{i}f_{i}h_i\\
        &=\sum\limits_{i=1}^{m}(u\parallel x_{i;t})_Fh_i,
    \end{align*}
  where the last equality follows from Proposition \ref{the relation of F-invariant and f-vector}.

\end{proof}
\begin{lemma}\label{F_invariant=0 for cluster variables}
    Let $\{z_1,\cdots,z_p\}$ be a set of unfrozen cluster variables of $\mathcal{A}$, where $p$ is a positive integer. If $(z_{i}\parallel z_{j})_F=0$ for any $i,j\in[1,p]$, then $\{z_1,\cdots,z_p\}$ is a subset of some cluster of $\mathcal{A}$.
\end{lemma}
\begin{proof}
    By Proposition \ref{the relation of F-invariant and f-vector} and Theorem \ref{relation of f-compatible and cluster}, there is a cluster containing $z_i$ and $z_j$ for any $i,j\in [1,p]$. Then the result follows from \cite[Theorem 4.3]{Cao_Li21} .
    \end{proof}
    As a consequence of Proposition \ref{liner F-invariant} and Lemma \ref{F_invariant=0 for cluster variables}, we have the following analogue of \cite[Theorem 4.22]{cao23} for generalized cluster algebras.
    
    \begin{proposition}
        \label{uu' is cluster monomial=(u,u')_F=0}
        Given two cluster monomials $u$ and $u'$. Then the product $u u'$ is still a cluster monomial in $\mathcal{U}$ if and only if $(u\parallel u')_F=0$.
        \end{proposition}
\bibliographystyle{plain} 
\bibliography{f-vector}
\end{document}